\documentclass[journal]{IEEEtran}
\usepackage[utf8]{inputenc}
\usepackage{amsmath}
\usepackage{amssymb}
\usepackage{bbm}
\usepackage{float}
\usepackage{setspace}
\usepackage{graphicx}
\usepackage{caption}
\usepackage{color}
\usepackage{pgfplots}
\usepackage{subdepth}
\usepackage{tikz}
\usetikzlibrary{arrows,positioning,decorations.pathreplacing,shapes,arrows,positioning,backgrounds}
\usepackage{algpseudocode,algorithm,algorithmicx}

\algrenewcommand\algorithmicrequire{\textbf{Input:}}
\algrenewcommand\algorithmicensure{\textbf{Output:}}
\newcommand{\norm}[1]{\left\lVert#1\right\rVert}
\newcommand*{\D}{\makebox[1ex]{\textbf{$\cdot$}}}%
\newtheorem{theorem}{Theorem}   

\newtheorem{problem}[theorem]{Problem}
\newtheorem{proposition}[theorem]{Proposition}
\newtheorem{corollary}[theorem]{Corollary}

\newtheorem{remark}{Remark}

\newtheorem{dfn}{Definition}
\DeclareMathOperator{\rank}{rank}
\DeclareMathOperator{\im}{im}

\DeclareMathOperator{\col}{col}

\tikzstyle{bball} = [circle,shading=ball, ball color=black!100!white,
minimum size=0.4cm]
\tikzstyle{wball} = [circle,shading=ball, ball color=white!100!black,
minimum size=0.4cm]
\begin{document} 
\title{Topology Reconstruction of Dynamical Networks via Constrained Lyapunov Equations}
\author{Henk~J.~van Waarde, Pietro~Tesi, and M.~Kanat~Camlibel%
\thanks{Henk van Waarde and Kanat Camlibel are with the Bernoulli Institute for Mathematics, Computer Science and Artificial Intelligence, Faculty of Science and Engineering, University of Groningen, P.O. Box 407, 9700 AK Groningen, The Netherlands. Pietro Tesi is with the Engineering and Technology Institute Groningen, Faculty of Science and Engineering, University of Groningen, 9747 AG Groningen, The Netherlands. Pietro Tesi is also with the Department of Information Engineering, University of Florence, 50139 Florence, Italy. Email: {\tt\small h.j.van.waarde@rug.nl, p.tesi@rug.nl, pietro.tesi@unifi.it, m.k.camlibel@rug.nl}}}%


\maketitle

\begin{abstract}
 The network structure (or topology) of a dynamical network is often unavailable or uncertain. Hence, we consider the problem of network reconstruction. Network reconstruction aims at inferring the topology of a dynamical network using measurements obtained from the network. In this technical note we define the notion of solvability of the network reconstruction problem. Subsequently, we provide necessary and sufficient conditions under which the network reconstruction problem is solvable. Finally, using constrained Lyapunov equations, we establish novel network reconstruction algorithms, applicable to general dynamical networks. We also provide specialized algorithms for specific network dynamics, such as the well-known consensus and adjacency dynamics.
\end{abstract}


\IEEEpeerreviewmaketitle

\section{Introduction}

\IEEEPARstart{N}{etworks} of dynamical systems appear in many contexts, including biological networks \cite{Sosa2005}, water distribution networks \cite{deRoo2015} and (wireless) sensor networks \cite{Mao2007}.

The overall behavior of a dynamical network is greatly influenced by its network structure (also called network topology). For instance, in the case of consensus networks, the dynamical network reaches \emph{consensus} if and only if the network graph is connected \cite{Saber2007}. Unfortunately, the interconnection structure of dynamical networks is often unavailable. For instance, in the case of wireless sensor networks \cite{Mao2007} the locations of sensors, and hence, communication links between sensors is not always known. Other examples of dynamical networks with unknown network topologies are encountered in biology, for instance in neural networks \cite{Sosa2005} and genetic networks \cite{Julius2009}. 

Consequently, the problem of \emph{network reconstruction} is studied in the literature. The aim of network reconstruction (also called topology identification) is to find the network structure and weights of a dynamical network, using measurements obtained from the network. To this end, most papers assume that the states of the network nodes can be measured. The literature on network reconstruction methods can roughly be divided into two parts, namely methods for \emph{stochastic} and \emph{deterministic} dynamical networks. 

Methods for stochastic network dynamics include \emph{inverse covariance estimation} \cite{Hassan-Moghaddam2016}, \cite{Morbidi2014} and methods based on \emph{power spectral analysis} \cite{Shahrampour2015}. Moreover, network reconstruction based on \emph{compressive sensing} \cite{Sanandaji2011} has been investigated. Furthermore, the authors of \cite{Materassi2012} consider network reconstruction using \emph{Wiener filtering}.

Apart from methods for stochastic networks, network reconstruction for deterministic network dynamics has been considered. In the paper \cite{Nabi-Abdolyousefi2010} the concept of \emph{node-knockout} is introduced, and a network reconstruction method based on this concept is discussed. The paper \cite{Goncalves2008} considers the problem of reconstructing a network topology from a transfer matrix of the network. Conditions are investigated under which the network structure can be uniquely determined. Furthermore, the paper \cite{Wu2016} considers network reconstruction using a so-called \emph{response network}.

In this note, we consider network reconstruction for \emph{deterministic} networks of linear dynamical systems. In contrast to papers studying network reconstruction for specific network dynamics such as consensus dynamics \cite{Nabi-Abdolyousefi2010} and adjacency dynamics \cite{Fazlyab2014}, we consider network reconstruction for \emph{general} linear network dynamics described by state matrices contained in the so-called \emph{qualitative class} \cite{Hogben2010}. It is our aim to infer the unknown network topology of such dynamical networks, from \emph{state measurements} obtained from the network.  

The contributions of this technical note are threefold. Firstly, we rigorously define what we mean by \emph{solvability} of the network reconstruction problem for dynamical networks. Loosely speaking, we say that the network reconstruction problem is solvable if the measurements obtained from a network correspond only with the network under consideration (and not with any other dynamical network). Secondly, we provide \emph{necessary} and \emph{sufficient} conditions under which the network reconstruction problem is solvable. Thirdly, we provide a framework for network reconstruction of dynamical networks, using constrained Lyapunov equations. We will show that our framework can be used to establish algorithms to infer network topologies for a variety of network dynamics, including Laplacian and adjacency dynamics. An attractive feature of our approach is that the conditions under which our algorithms reconstruct the network structure are not restrictive. In other words, we show that our algorithms return the correct network structure if and only if the network reconstruction problem is solvable. 

Although this note mainly focuses on continuous-time network dynamics, we also show how our reconstruction algorithms can be applied to discrete-time systems, and to systems with sampled measurements. 

The organization of this technical note is as follows. First, in Section \ref{section preliminaries}, we introduce preliminaries and notation used in this note. Subsequently, we give a formal problem statement in Section \ref{section problem formulation}. In Section \ref{section results solvability} we discuss necessary and sufficient conditions for the solvability of the network reconstruction problem. Section \ref{section results identification} provides our network reconstruction algorithms. We consider an illustrative example in Section \ref{section illustrative example}. Finally, Section \ref{section conclusions} contains our conclusions.

\section{Preliminaries}
\label{section preliminaries}
We denote the set of natural, real, and complex numbers by $\mathbb{N}$, $\mathbb{R}$, and $\mathbb{C}$ respectively. Moreover, the set of real $m \times n$ matrices is denoted by $\mathbb{R}^{m \times n}$. We denote the set of positive (non-negative) real numbers by $\mathbb{R}_{> 0}$ (respectively, $\mathbb{R}_{\geq 0}$). Furthermore, the set of all symmetric $n \times n$ matrices is given by $\mathbb{S}^n$. The vector of ones is denoted by $\mathbbm{1}$. Furthermore, for $x_1,x_2,\dots,x_n \in \mathbb{R}$, we use the notation $\col(x_1,x_2,\dots,x_n) \in \mathbb{R}^n$, which denotes the $n$-dimensional column vector with elements $x_1,x_2,\dots,x_n$. The \emph{image} of a matrix $A$ is denoted by $\im A$ and the \emph{kernel} of $A$ is denoted by $\ker A$. For a given set $S$, the \emph{power set} $2^S$ is the set of all subsets of $S$. Let $X$ and $Y$ be nonempty sets. If for each $x \in X$, there exists a set $F(x) \subseteq Y$, we say $F$ is a \emph{set-valued map} from $X$ to $Y$, and we denote $F : X \to 2^Y$. The \emph{image} of a set-valued map $F: X \to 2^Y$ is defined as $\im F := \{ y \in Y \mid \exists x \in X \text{ such that } y \in F(x) \}$.

\subsection{Preliminaries on systems theory}
\label{preliminaries systems}
Consider the linear time-invariant system
\begin{equation}
\label{linear system}
\begin{aligned}
\dot{x}(t) &= Ax(t) \\
y(t) &= Cx(t),
\end{aligned}
\end{equation}
where $x \in \mathbb{R}^n$ is the state, $y \in \mathbb{R}^p$ is the output, and the real matrices $A$ and $C$ are of suitable dimensions. We denote the \emph{unobservable subspace} of system \eqref{linear system} by $\big \langle \ker C \: | \: A \big \rangle$, i.e.,
\begin{equation*}
\big \langle \ker C \: | \: A \big \rangle := \bigcap\limits_{i=0}^{n-1} \ker \left( CA^i \right).
\end{equation*}
The subspace $\big \langle \ker C \: | \: A \big \rangle$ is $A$-invariant, that is, $A \big \langle \ker C \: | \: A \big \rangle \subseteq \big \langle \ker C \: | \: A \big \rangle$. Furthermore, system \eqref{linear system} is \emph{observable} if and only if $\big \langle \ker C \: | \: A \big \rangle = \{ 0 \}$ (see, e.g., Chapter 3 of \cite{Trentelman2001}). If system \eqref{linear system} is observable, we say the pair $(C,A)$ is observable. 
    
\subsection{Preliminaries on graph theory}
\label{preliminaries graphs}
All graphs considered in this note are simple, i.e., without self-loops and with at most one edge between any pair of vertices. We denote the set of simple, undirected graphs of $n$ nodes by $\mathcal{G}^n$. Consider a graph $G \in \mathcal{G}^n$, with vertex set $V = \{1,2,\dots,n\}$ and edge set $E$. The set of neighbors $\mathcal{N}_i$ of vertex $i \in V$ is defined as $\mathcal{N}_i := \{j \in V \mid (i,j) \in E \}$. 

We will now define various families of matrices associated with graphs in $\mathcal{G}^n$. To this end, we first define the set-valued map $\mathcal{Q} : \mathcal{G}^n \to 2^{\mathbb{S}^n}$ as 
\begin{equation*}
\mathcal{Q}(G) := \{X \in \mathbb{S}^n \mid \text{for all } i \neq j, \: X_{ij} \neq 0 \iff (i,j) \in E \}.
\end{equation*}
The set of matrices $\mathcal{Q}(G)$ is called the \emph{qualitative class} of the graph $G \in \mathcal{G}^n$ \cite{Hogben2010}. The qualitative class has recently been studied in the context of structural controllability of dynamical networks \cite{Monshizadeh2014}, \cite{vanWaarde2017}. Note that each matrix $X \in \mathcal{Q}(G)$ carries the graph structure of $G$, in the sense that $X$ contains nonzero off-diagonal entries in exactly the same positions corresponding to the edges in $G$. Furthermore, note that the diagonal elements of matrices in the qualitative class are unrestricted. Hence, examples of matrices in $\mathcal{Q}(G)$ include the well-known (weighted) \emph{adjacency} and \emph{Laplacian} matrices, which are defined next. Define the set-valued map $\mathcal{A} : \mathcal{G}^n \to 2^{\mathbb{S}^n}$ as 
\begin{equation*}
\mathcal{A}(G) := \{A \in \mathcal{Q}(G) \mid \text{for all } i,j, \: A_{ij} \geq 0 \text{ and } A_{ii} = 0 \}.
\end{equation*}
Matrices in $\mathcal{A}(G)$ are called adjacency matrices associated with the graph $G$. Subsequently, define the set-valued map $\mathcal{L} : \mathcal{G}^n \to 2^{\mathbb{S}^n}$ as  
\begin{equation*}
\mathcal{L}(G) := \{L \in \mathcal{Q}(G) \mid L \mathbbm{1} = 0 \text{ and for all } i \neq j, L_{ij} \leq 0 \}.
\end{equation*}
Matrices in the set $\mathcal{L}(G)$ are called Laplacian matrices of $G$. A Laplacian matrix $L \in \mathcal{L}(G)$ is said to be \emph{unweighted} if $L_{ij} \in \{0,-1\}$ for all $i \neq j$. Similarly, an adjacency matrix $A \in \mathcal{A}(G)$ is called unweighted if $A_{ij} \in \{0,1\}$ for all $i, j$.

\subsection{Preliminaries on consensus dynamics}
\label{preliminaries consensus}
Consider a graph $G \in \mathcal{G}^n$, with vertex set $V = \{1,2,\dots,n\}$ and edge set $E$. With each vertex $i \in V$, we associate a linear dynamical system $
\dot{x}_i(t) = u_i(t)$, where $x_i \in \mathbb{R}$ is the state of node $i$, and $u_i \in \mathbb{R}$ is its control input. Suppose that each node $i \in V$ applies the control input
\begin{equation*}
u_i(t) = - \sum_{j \in \mathcal{N}_i} a_{ij} (x_i(t) - x_j(t)),
\end{equation*}
where $a_{ij} = a_{ji} > 0$ for all $i \in V$ and $j \in \mathcal{N}_i$. Then, the dynamics of the overall system can be written as 
\begin{equation}
\label{consensus}
\dot{x}(t) = -L x(t),
\end{equation}
where $x = \col(x_1,x_2,\dots,x_n)$, and $L \in \mathcal{L}(G)$ is a Laplacian matrix. We refer to system \eqref{consensus} as a \emph{consensus network}. Consensus networks have been studied extensively in the literature, see, e.g., \cite{Saber2007} and the references therein.

\section{Problem formulation}
\label{section problem formulation}
In this section we define the network reconstruction problem. We consider a linear time-invariant network system, with nodes satisfying single-integrator dynamics. We assume that the state matrix of the system (and hence, the network topology) is not directly available. Moreover, we suppose that the state vector of the system is available for measurement during a time interval $[0,T]$. It is our goal to find conditions on the system under which the exact state matrix can be reconstructed from such measurements. Moreover, if the state matrix can be reconstructed, we want to develop algorithms to infer the state matrix from measurements. 

We will now make these problems more precise. Since we want to consider network reconstruction for general network dynamics (instead of specific consensus or adjacency dyna\-mics), we consider any set-valued map $\mathcal{K} : \mathcal{G}^n \to 2^{\mathbb{S}^n}$ such that for all $G \in \mathcal{G}^n$, we have that
\begin{equation}
\label{mappingK}
\mathcal{K}(G) \subseteq \mathcal{Q}(G)
\end{equation}
is a nonempty subset. The map $\mathcal{K}$ is specified by the available information on the \emph{type} of network. For example, if we know that we deal with a consensus network, we have $\mathcal{K} = -\mathcal{L}$. On the other hand, if no additional information on the communication weights (such as sign constraints) is known, we let $\mathcal{K} = \mathcal{Q}$. With this in mind, we consider the system
\begin{equation}
\begin{aligned}
\label{networked system}
\dot{x}(t) &= X x(t) \text{ for } t \in \mathbb{R}_{\geq 0} \\
x(0) &= x_0,
\end{aligned}
\end{equation}
where $x \in \mathbb{R}^n$ is the state, and $X \in \im \mathcal{K}$ (i.e., $X \in \mathcal{K}(G)$ for some network graph $G \in \mathcal{G}^n$). In what follows, we denote the state trajectory of \eqref{networked system} by $x_{x_0}(\D)$, where the subscript indicates dependence on the initial condition $x_0$. We assume that $X$ is unknown, but the state trajectory of \eqref{networked system} can be measured during the time-interval $[0,T]$, where $T \in \mathbb{R}_{> 0}$. The problem of \emph{network reconstruction} concerns finding the matrix $X$ (and thereby, the graph $G$), using the state measurements $x_{x_0}(t)$ for $t \in [0,T]$. Of course, this is only possible if the state trajectory $x_{x_0}(\D)$ of \eqref{networked system} is not a solution to the differential equation $\dot{x}(t) = \bar{X} x(t)$ for some other admissible state matrix $\bar{X} \neq X$. Indeed, if this were the case, the state measurements could correspond to a network described either by $X$ or $\bar{X}$, and we would not be able to distinguish between the two. This leads to the following definition. 
\begin{dfn}
	\label{definition solvability}	
	Consider system \eqref{networked system}, and denote its state trajectory by $x_{x_0}(\D)$. We say that the network reconstruction problem is \emph{solvable} for system \eqref{networked system} if for all $\bar{X} \in \im \mathcal{K}$ such that $x_{x_0}(\D)$ is a solution to 
	\begin{equation}
	\label{other admissible}
	\dot{x}(t) = \bar{X} x(t) \text{ for } t \in [0,T],
	\end{equation}
	we have $\bar{X} = X$. In the case that the network reconstruction problem is solvable for system \eqref{networked system}, we say that the network reconstruction problem is solvable for $(x_0,X,\mathcal{K})$.
\end{dfn}

\begin{remark}
	\label{remark analytic}
	As the state variables of system \eqref{networked system} are sums of exponential functions of $t$, they are \emph{real analytic} functions of $t$. It is well-known that if two real analytic functions are equal on a non-degenerate interval, they are equal on their whole domain (see, e.g., Corollary 1.2.5 of \cite{Krantz2002}). Consequently, the state vector $x_{x_0}(\D)$ of system \eqref{networked system} satisfies \eqref{other admissible} for $t \in [0,T]$ if and only if $x_{x_0}(\D)$ satisfies \eqref{other admissible} for $t \in \mathbb{R}_{\geq 0}$. Therefore, Definition \ref{definition solvability} can be equivalently stated for $t \in \mathbb{R}_{\geq 0}$ instead of $t \in [0,T]$. 
\end{remark}

In this note we are interested in conditions on $x_0$, $X$, and $\mathcal{K}$ under which the network reconstruction problem is solvable for $(x_0,X,\mathcal{K})$. More explicitly, we have the following problem. 
\begin{problem}
	\label{problem conditions}
	Consider system \eqref{networked system}. Provide necessary and sufficient conditions on $x_0$, $X$, and $\mathcal{K}$ under which the network reconstruction problem is solvable for system \eqref{networked system}. 
\end{problem}
In addition to Problem \ref{problem conditions}, we are interested in solving the network reconstruction problem itself. This is stated in the following problem.
\begin{problem}
	\label{problem identify}
	Consider system \eqref{networked system}, and denote its state vector by $x_{x_0}(\D)$. Suppose that $x_{x_0}(\D)$ is available for measurement during the time interval $[0,T]$, and that the network reconstruction problem is solvable for \eqref{networked system}. Provide a method to compute the matrix $X$. 
\end{problem}

\begin{remark}
	Note that we assume that the states of all nodes in the network can be measured. This assumption is \emph{necessary} in the sense that the network reconstruction problem is not solvable (in the case of $\mathcal{Q}(G)$) if we can only measure a part of the state vector. To see this, suppose that we only have access to a $p$-dimensional output vector $y(t) = C x(t)$, where
	\begin{equation*}
	C = \begin{pmatrix}
	I & 0
	\end{pmatrix} \in \mathbb{R}^{p \times n}.
	\end{equation*}
	We claim that for each $X \in \mathcal{Q}(G)$ and $x_0 \in \mathbb{R}^n$ there exists a graph $\bar{G}$, a matrix $\bar{X} \in \mathcal{Q}(\bar{G}) \setminus \{X\}$ and a vector $\bar{x}_0 \in \mathbb{R}^n$ such that 
	\begin{equation}
	\label{outputequal}
	Ce^{Xt}x_0 = Ce^{\bar{X}t}\bar{x}_0.
	\end{equation}
	That is, we cannot distinguish between $X$ and $\bar{X}$ on the basis of output measurements. To see that this claim is true, we write $X$ as
	\begin{equation*}
	X = \begin{pmatrix}
	X_{11} & X_{12} \\ X_{21} & X_{22} 
	\end{pmatrix},
	\end{equation*}
	where the partitioning of $X$ is compatible with the one of $C$. Now we distinguish two cases. First suppose that $X_{21} \neq 0$. Clearly, there exists a vector $z \in \mathbb{R}^{n-p}$ such that $z^\top X_{21} \neq 0$ and $z^\top z = 1$. Define 
	\begin{equation*}
	S := \begin{pmatrix}
	I & 0 \\ 0 & S_{22}
	\end{pmatrix},
	\end{equation*}
	where $S_{22} := I - 2 z z^\top = S_{22}^{-1}$. Then let $\bar{X} := S X S$ and $\bar{x}_0 := S x_0$. It is not difficult to see that $\bar{X} \neq X$ and \eqref{outputequal} is satisfied for this choice of $\bar{X}$ and $\bar{x}_0$. Secondly, consider the case that $X_{21} = 0$. Then $X_{12} = 0$. We can choose $\bar{X}$ as
	\begin{equation*}
	\bar{X} := \begin{pmatrix}
	X_{11} & 0 \\ 0 & \bar{X}_{22}
	\end{pmatrix},
	\end{equation*}
	where $\bar{X}_{22} \neq X_{22}$. In this case, it can be shown that $\bar{X}$ and $\bar{x}_0 := x_0$ satisfy \eqref{outputequal}. Hence, for the network reconstruction problem to be solvable it is necessary to measure all nodes.
\end{remark}

\section{Main results: solvability of the network reconstruction problem}
\label{section results solvability}
In this section we state our main results regarding Problem \ref{problem conditions}. That is, we provide conditions on $x_0$, $X$, and $\mathcal{K}$ under which the network reconstruction problem is solvable. Firstly, in Section \ref{section solvability K} we provide necessary and sufficient conditions for the solvability of the network reconstruction problem in the general case that $\mathcal{K}$ is any mapping satisfying \eqref{mappingK}. Later, we consider the special cases in which $\mathcal{K} = \mathcal{Q}$ (Section \ref{section solvability Q}), and the cases in which $\mathcal{K} = -\mathcal{L}$ or $\mathcal{K} = \mathcal{A}$ (Section \ref{section solvability L}).

\subsection{Solvability for general $\mathcal{K}$}
\label{section solvability K}
In this section, we provide a general solution to Problem \ref{problem conditions}. Let $G \in \mathcal{G}^n$ be a graph, and let the mapping $\mathcal{K}$ be as in \eqref{mappingK}. Recall that we consider the dynamical network described by system \eqref{networked system}. As a preliminary result, we give conditions under which the state trajectory $x_{x_0}(\D)$ of system \eqref{networked system} is also the solution to the system 
\begin{equation}
\label{networked system2}
\begin{aligned}
\dot{x}(t) &= \bar{X} x(t) \text{ for } t \in \mathbb{R}_{\geq 0} \\
x(0) &= x_0,
\end{aligned}
\end{equation}
where $\bar{X} \in \mathcal{K}(\bar{G})$ for some graph $\bar{G} \in \mathcal{G}^n$. This result is given in the following proposition. 
\begin{proposition}
	\label{propositiontwomatrices}
	Consider systems \eqref{networked system} and \eqref{networked system2}, and let $x_{x_0}(\D)$ be the state trajectory of \eqref{networked system}. The trajectory $x_{x_0}(\D)$ is also the solution to system \eqref{networked system2} if and only if $x_0 \in \big \langle \ker\left( \bar{X} - X \right) \: | \: X \big \rangle$.
\end{proposition}

\begin{IEEEproof}
	Suppose that the state trajectory $x_{x_0}(\D)$ of \eqref{networked system} is also the solution to system \eqref{networked system2}. This means that $x_{x_0}(\D)$ is the solution to both the differential equation
	\begin{equation}
	\label{nominal system}
	\dot{x}(t) = X x(t) \text{ for } t \in \mathbb{R}_{\geq 0},
	\end{equation}
	and the differential equation
	\begin{equation}
	\label{networked system3}
	\dot{x}(t) = X x(t) + (\bar{X} - X) x(t) \text{ for } t \in \mathbb{R}_{\geq 0}.
	\end{equation}
	In particular, by substitution of $t = 0$, this implies that $x_0$ is contained in $\ker\left(\bar{X} - X\right)$. Moreover, by taking the $i$-th time-derivative of \eqref{nominal system} and \eqref{networked system3}, we find that $x_0 \in \ker\left(\bar{X} - X\right)X^i$ for $i = 1,2,\dots,n-1$. Consequently, we obtain $x_0 \in \big \langle \ker\left( \bar{X} - X \right) \: | \: X \big \rangle$. 
	
	\noindent
	Conversely, suppose that the initial state $x_0$ of system \eqref{networked system} satisfies $x_0 \in \big \langle \ker\left( \bar{X} - X \right) \: | \: X \big \rangle$. By $X$-invariance of $\big \langle \ker\left( \bar{X} - X \right) \: | \: X \big \rangle$, this implies that the state trajectory $x_{x_0}(\D)$ of system \eqref{networked system} satisfies $x_{x_0}(t) \in \big \langle \ker\left( \bar{X} - X \right) \: | \: X \big \rangle$ for all $t \in \mathbb{R}_{\geq 0}$. Specifically, we have that $x_{x_0}(t) \in \ker\left(\bar{X} - X\right)$ for all $t \in \mathbb{R}_{\geq 0}$. We conclude that $x_{x_0}(\D)$ is the solution to Equation \eqref{networked system3}, and consequently, to Equation \eqref{networked system2}.
\end{IEEEproof}
\begin{remark}
	Note that a condition equivalent to the one given in Proposition \ref{propositiontwomatrices} can be stated in terms of the \emph{common eigenspaces} of $X$ and $\bar{X}$. Such a condition was previously proven by Battistelli \emph{et al.} \cite{Battistelli2018}, \cite{Battistelli2015} in the case that $X$ and $\bar{X}$ are \emph{Laplacian} matrices.
\end{remark}
By combining Proposition \ref{propositiontwomatrices} and the fact that the state variables of \eqref{networked system} are real analytic functions in $t$ (see Remark \ref{remark analytic}), we obtain Theorem \ref{theorem nec suf}. This theorem states a necessary and sufficient condition under which the network reconstruction problem is solvable for $(x_0,X,\mathcal{K})$. 
\begin{theorem}
	\label{theorem nec suf}
	Let $G \in \mathcal{G}^n$ be a graph, and let the mapping $\mathcal{K}$ be as in \eqref{mappingK}. Moreover, consider a matrix $X \in \mathcal{K}(G)$ and a vector $x_0 \in \mathbb{R}^n$. The network reconstruction problem is solvable for $(x_0,X,\mathcal{K})$ if and only if for all $\bar{X} \in \im \mathcal{K} \setminus \{X\}$, we have
	\begin{equation*}
	x_0 \not\in \big \langle \ker\left( \bar{X} - X \right) \: | \: X \big \rangle.
	\end{equation*}
\end{theorem}
\begin{remark}
Although Theorem \ref{theorem nec suf} gives a general necessary and sufficient condition for network reconstruction, it is not directly clear how to verify this condition. Especially since $X$ is assumed to be unknown, it seems difficult to check that $x_0 \not\in \big \langle \ker\left( \bar{X} - X \right) \: | \: X \big \rangle$. In fact, we will show in Section \ref{section results identification} that the condition of Theorem \ref{theorem nec suf} can be checked using only the measurements $x_{x_0}(t)$ for $t \in [0,T]$.
\end{remark}
Note that the condition of Theorem \ref{theorem nec suf} is not only given in terms of $x_0$ and $X$, but also in terms of all other matrices $\bar{X} \in \im \mathcal{K}$. In the following theorem, we provide a simple \emph{sufficient} condition for the solvability of the network reconstruction problem, which is stated in terms of $x_0$ and $X$.

\begin{theorem}
	\label{sufficiency observability}
	Let $G \in \mathcal{G}^n$ be a graph, and let the mapping $\mathcal{K}$ be as in \eqref{mappingK}. Moreover, consider a matrix $X \in \mathcal{K}(G)$ and a vector $x_0 \in \mathbb{R}^n$. The network reconstruction problem is solvable for $(x_0,X,\mathcal{K})$ if the pair $(x_0^\top,X)$ is observable.
\end{theorem}
\begin{IEEEproof}
	Suppose that the pair $(x_0^\top,X)$ is observable, and assume that $x_0 \in \big \langle \ker\left( \bar{X} - X \right) \: | \: X \big \rangle$ for some matrix $\bar{X} \in \im \mathcal{K}$. We want to show that $\bar{X} = X$. Note that by hypothesis, we have $x_0 \in \ker(\bar{X} - X)X^i$, for $i = 0,1,\dots,n-1$. As a consequence, we obtain the equalities 
	\begin{equation}
	\label{n equalities}
	\bar{X}X^i x_0 = X^{i+1} x_0,
	\end{equation}
	for $i = 0,1,\dots,n-1$. It is not difficult to see that by induction, Equation \eqref{n equalities} implies that
	\begin{equation}
	\label{powers of X}
	X^i x_0 = \bar{X}^i x_0,
	\end{equation}
	for $i = 1,2,\dots,n$. In other words, the matrix $X \begin{pmatrix}
	x_0 & X x_0 & \dots & X^{n-1} x_0 
	\end{pmatrix}$ is equal to 
	\begin{equation}
	\label{powers of X2}
	\bar{X} \begin{pmatrix}
	x_0 & \bar{X} x_0 & \dots & \bar{X}^{n-1} x_0 
	\end{pmatrix}.
	\end{equation}
	Since the pair $(x_0^\top,X)$ is observable and $X$ is symmetric, the matrix $\begin{pmatrix}
	x_0 & X x_0 & \dots & X^{n-1} x_0 
	\end{pmatrix}$ is invertible. This allows us to conclude that $X$ equals
	\begin{equation*}
	\bar{X} \begin{pmatrix}
	x_0 & \bar{X} x_0 & \dots & \bar{X}^{n-1} x_0 
	\end{pmatrix}\begin{pmatrix}
	x_0 & X x_0 & \dots & X^{n-1} x_0 
	\end{pmatrix}^{-1}.
	\end{equation*}
	However, by \eqref{powers of X}, this implies that $X = \bar{X}$. Consequently, for all $\bar{X} \in \im \mathcal{K} \setminus \{X\}$ we have $x_0 \not\in \big \langle \ker\left( \bar{X} - X \right) \: | \: X \big \rangle$. Finally, we conclude by Theorem \ref{theorem nec suf} that the network reconstruction problem is solvable for $(x_0,X,\mathcal{K})$.
\end{IEEEproof}
In the next section, we show that for a specific mapping $\mathcal{K}$, the observability condition of Theorem \ref{sufficiency observability} is necessary \emph{and} sufficient. However, in general, the observability condition is not necessary. In particular, this will be shown for consensus networks in Section \ref{section solvability L}.

\subsection{Solvability for $\mathcal{K} = \mathcal{Q}$}
\label{section solvability Q}
In this subsection, we consider the case that $\mathcal{K} = \mathcal{Q}$. This case corresponds to the situation where we do not have any additional information (such as sign constraints) on the entries of the state matrix $X$. To be precise, we consider system \eqref{networked system},
where $X \in \mathcal{Q}(G)$ for some network graph $G \in \mathcal{G}^n$. We will see that the solvability of the network reconstruction problem for $(x_0,X,\mathcal{Q})$ is in fact equivalent to the observability of the pair $(x_0^\top,X)$. This is stated in the following theorem. 
\begin{theorem}
	\label{theorem Q nec suf}
	Consider a graph $G \in \mathcal{G}^n$, let $X \in \mathcal{Q}(G)$, and let $x_0 \in \mathbb{R}^n$. The network reconstruction problem is solvable for $(x_0,X,\mathcal{Q})$ if and only if the pair $(x_0^\top,X)$ is observable.
\end{theorem}

\begin{IEEEproof}
	Sufficiency follows immediately from Theorem \ref{sufficiency observability} by taking $\mathcal{K} = \mathcal{Q}$. Hence, assume that the pair $(x_0^\top,X)$ is unobservable. We want to show that the network reconstruction problem is not solvable for $(x_0,X,\mathcal{Q})$. To do so, we will construct a matrix $\bar{X} \neq X$ such that $x_0 \in \big \langle \ker\left( \bar{X} - X \right) \: | \: X \big \rangle$. 
	
	\noindent
	Let $v \in \mathbb{R}^n$ be a nonzero vector such that 
	\begin{equation}
	v^\top \begin{pmatrix}
	x_0 & X x_0 & \dots & X^{n-1}x_0 \end{pmatrix} = 0.
	\end{equation}
	Such a vector exists, as $(x_0^\top,X)$ is unobservable. Subsequently, define the matrix $\bar{X} := X + v v^\top$. By definition of $v$, we obtain $\bar{X}^i x_0 = X^i x_0$, for $i = 1,2,\dots,n$. Consequently, $x_0 \in \big \langle \ker\left( \bar{X} - X \right) \: | \: X \big \rangle$. It remains to be shown that $\bar{X} \in \im \mathcal{Q}$, i.e., $\bar{X} \in \mathcal{Q}(\bar{G})$ for some $\bar{G} \in \mathcal{G}^n$. Define the simple undirected graph $\bar{G} = (V,E)$, where $V := \{1,2,\dots,n\}$, and for distinct $i,j \in V$, we have $(i,j) \in E$ if and only if $\bar{X}_{ij} \neq 0$. By definition of the qualitative class $\mathcal{Q}(\bar{G})$, we obtain $\bar{X} \in \mathcal{Q}(\bar{G})$. We conclude that the network reconstruction problem is not solvable for $(x_0,X,\mathcal{Q})$.
\end{IEEEproof}

\subsection{Solvability for $\mathcal{K} = -\mathcal{L}$ and $\mathcal{K} = \mathcal{A}$}
\label{section solvability L}
In what follows, we consider solvability of the network reconstruction problem for consensus and adjacency networks. We will start with consensus networks. That is, we consider the system
\begin{equation}
\begin{aligned}
\label{networked system5}
\dot{x}(t) &= -L x(t) \text{ for } t \in \mathbb{R}_{\geq 0} \\
x(0) &= x_0,
\end{aligned}
\end{equation}
where $x \in \mathbb{R}^n$ is the state and $L \in \mathcal{L}(G)$ denotes the Laplacian matrix of a graph $G \in \mathcal{G}^n$. In this section we show by means of an example that observa\-bility of $(x_0^\top,-L)$ is not necessary for the solvability of the network reconstruction problem for $(x_0,-L,-\mathcal{L})$. In Section \ref{section results identification} we will use this fact to establish an algorithm for network reconstruction of consensus networks, that does not require observability of the pair $(x_0^\top,-L)$. Consider the star graph $G$ and Laplacian matrix $L$, depicted in Figure \ref{fig:star graph}.

\noindent\begin{minipage}{.48\linewidth}
\centering
\includegraphics[height=3cm]{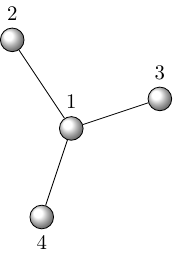}
\end{minipage}%
\begin{minipage}{.48\linewidth}
	\begin{equation*}
	L = \begin{pmatrix}
	\phantom{-}3 & -1 & -1 & -1 \\
	-1 & \phantom{-}1 & \phantom{-}0 & \phantom{-}0 \\
	-1 & \phantom{-}0 & \phantom{-}1 & \phantom{-}0 \\
	-1 & \phantom{-}0 & \phantom{-}0 & \phantom{-}1
	\end{pmatrix}.
	\end{equation*}
\end{minipage}
\captionof{figure}{Star graph $G$ with Laplacian matrix $L$.}
\label{fig:star graph} 
\vskip.2\baselineskip

Moreover, consider the initial condition $x_0 \in \mathbb{R}^4$ given by $x_0 = \col(1,0,3,1)$. We claim that the network reconstruction problem is solvable for $(x_0,-L,-\mathcal{L})$, even though the pair $(x_0^\top,-L)$ is unobser\-vable. Indeed, it can be verified that the unobservable subspace of $(x_0^\top,-L)$ is $\big \langle \ker x_0^\top \: | \: -L \big \rangle = \im v$,
where the vector $v$ is defined as $v := \col(0,2,1,-3)$. This implies that $(x_0^\top,-L)$ is unobservable. To prove that the network reconstruction problem is solvable for $(x_0,-L,-\mathcal{L})$, consider a Laplacian matrix $\bar{L} \in \im \mathcal{L}$ such that $x_0 \in \big \langle \ker\left( L - \bar{L} \right) \: | \: -L \big \rangle$. Following the proof of Theorem \ref{sufficiency observability} (Equations \eqref{powers of X} and \eqref{powers of X2}), we find that 
\begin{equation}
(L-\bar{L}) \begin{pmatrix}
x_0 & -L x_0 & \dots & (-L)^{n-1} x_0 
\end{pmatrix}
= 0.
\end{equation}
In other words, the columns of the matrix $D := L-\bar{L}$ are contained in the unobservable subspace of $(x_0^\top,-L)$. Since $D$ is symmetric and $\big \langle \ker x_0^\top \: | \: -L \big \rangle = \im v$, we find
\begin{equation}
D = \alpha v v^\top = \alpha \begin{pmatrix}
0 & \phantom{-}0 & \phantom{-}0 & \phantom{-}0 \\
0 & \phantom{-}4 & \phantom{-}2 & -6 \\
0 & \phantom{-}2 & \phantom{-}1 & -3 \\
0 & -6 & -3 & \phantom{-}9
\end{pmatrix},
\end{equation} 
for some $\alpha \in \mathbb{R}$. If $\alpha \neq 0$, the entries $D_{32}$ and $D_{42}$ of the matrix $D$ have opposite sign. Since we have $L_{32} = L_{42} = 0$, we conclude from the relation $\bar{L} = L - D$ that $\bar{L}_{32}$ and $\bar{L}_{42}$ have opposite sign. However, this is a contradiction as $\bar{L}$ is a Laplacian matrix. Therefore, we conclude that $\alpha = 0$, and hence, $D = 0$. Consequently, we obtain $L = \bar{L}$. By Theorem \ref{theorem nec suf}, we conclude that the network reconstruction problem is solvable for $(x_0,-L,-\mathcal{L})$. Thus, we have shown that obser\-vability of the pair $(x_0^\top,L)$ is not necessary for the solvability of the network reconstruction problem for $(x_0,-L,-\mathcal{L})$.

It can be shown that $\im v$ also equals $\big \langle \ker x_0^\top \: | \: A \big \rangle$, where $A \in \mathcal{A}(G)$ denotes the unweighted \emph{adjacency} matrix asso\-ciated with the star graph $G$ depicted in Figure \ref{fig:star graph}. Then, using the exact same reasoning as before, we conclude that the pair $(x_0^\top,A)$ is unobservable, but the network reconstruction problem is solvable for $(x_0,A,\mathcal{A})$. In other words, observability of $(x_0^\top,A)$ is not necessary for the solvability of the network reconstruction problem for $(x_0,A,\mathcal{A})$.

\section{Main results: the network reconstruction problem}
\label{section results identification}

In this section, we provide a solution to Problem \ref{problem identify}. That is, given measurements generated by an unknown network, we establish algorithms to infer the network topology. Similar to the setup of Section \ref{section results solvability}, we start with the most general case in which  $\mathcal{K}$ is any mapping satisfying \eqref{mappingK}. For this case, we obtain a general methodology to infer $X \in \mathcal{K}(G)$ from measurements. Subsequently, we provide specific algorithms for network reconstruction in the case that $\mathcal{K} = \mathcal{Q}$ (Section \ref{section identification Q}), and in the case of consensus and adjacency networks (Section \ref{section identification L}). 

\subsection{Network reconstruction for general $\mathcal{K}$}
\label{section identification K}
Recall that we consider the system \eqref{networked system}, where the matrix $X$ and graph $G$ are unknown, but the state vector $x_{x_0}(\D)$ of \eqref{networked system} can be measured during the time interval $[0,T]$. In this section, we establish a method to infer the matrix $X$ and graph $G$ using the vector $x_{x_0}(t)$ for $t \in [0,T]$. Firstly, define the matrix
\begin{equation}
\label{transformed P}
P := \int_0^{T} x_{x_0}(t) x_{x_0}(t)^\top dt = \int_0^{T} e^{Xt} x_0 x_0^\top e^{Xt} dt.
\end{equation}
Note that $P$ can be computed from the measurements $x_{x_0}(t)$ for $t \in [0,T]$. The unknown matrix $X$ is a solution to a \emph{Lyapunov equation} involving the matrix $P$. Indeed, we have 
\begin{equation}
\label{derivation}
\begin{aligned}
XP + PX &= \int_0^{T} \left( X e^{Xt} x_0 x_0^\top e^{Xt} + e^{Xt} x_0 x_0^\top e^{Xt} X \right) dt \\
&= \int_0^{T} \frac{d}{dt} \left( e^{Xt} x_0 x_0^\top e^{Xt} \right) dt \\
&= x_T x_T^\top - x_0 x_0^\top,
\end{aligned}
\end{equation}
where $x_T := x_{x_0}(T) = e^{X T} x_0$. In other words, $X$ satisfies the Lyapunov equation 
\begin{equation}
\label{lyapunov equation transformed}
\begin{aligned}
XP + PX &= Q,
\end{aligned}
\end{equation}
where $Q$ is defined as $Q := x_T x_T^\top - x_0 x_0^\top$. Note that we can compute the matrix $Q$ from the measurements $x_{x_0}(t)$ at time $t = 0$ and time $t = T$. Therefore, if the matrix $S = X$ is the \emph{unique} solution to the Lyapunov equation 
\begin{equation}
\label{lyapunov equation transformed2}
\begin{aligned}
SP + PS &= Q,
\end{aligned}
\end{equation}
we can find $X$ (and therefore $G$), by solving \eqref{lyapunov equation transformed2} for $S$. However, it turns out that in general it is not necessary for network reconstruction that the Lyapunov equation \eqref{lyapunov equation transformed2} has a unique solution $S$. In fact, we only need a unique solution $S$ in the image of $\mathcal{K}$. That is, the Lyapunov equation \eqref{lyapunov equation transformed2} may have many solutions, but if only one of these solutions is contained in $\im \mathcal{K}$, we can solve the network reconstruction problem for $(x_0,X,\mathcal{K})$. This is stated more formally in the following theorem. 

\begin{theorem}
	\label{theorem unique solution}
	Let $G \in \mathcal{G}^n$, and let the mapping $\mathcal{K}$ be as in \eqref{mappingK}. Moreover, consider $X \in \mathcal{K}(G)$, $x_0 \in \mathbb{R}^{n}$, and let $P$ and $Q$ be as defined in \eqref{transformed P} and \eqref{lyapunov equation transformed} respectively. The network reconstruction problem is solvable for $(x_0,X,\mathcal{K})$ if and only if there exists a unique matrix $S$ satisfying 
	\begin{equation}
	\label{lyapunov equation theorem1}
	SP + PS = Q, \quad S \in \im \mathcal{K}.
	\end{equation}
	Moreover, under this condition, we have $S = X$. 
\end{theorem}
Before we can prove Theorem \ref{theorem unique solution}, we need the following proposition, which states that $\ker P$ equals the unobservable subspace of the pair $(x_0^\top,X)$.

\begin{proposition}
	\label{proposition kernel}
	Let $P$, $x_0$ and $X$ be as in \eqref{transformed P}. Then we have $\ker P = \big \langle \ker x_0^\top \: | \: X \big \rangle$.
\end{proposition}

\begin{IEEEproof}
Let $v \in \ker P$. We have 
\begin{equation}
v^\top P v = \int_0^\top \left( x_0^\top e^{Xt} v \right)^2 dt = 0, 
\end{equation}
from which we obtain $x_0^\top e^{Xt} v = 0$ for all $t \in [0,T]$. Since $x_0^\top e^{Xt} v$ is a real analytic function, we see that $x_0^\top e^{Xt} v = 0$ for all $t \in \mathbb{R}_{\geq 0}$ (cf. Remark \ref{remark analytic}). This implies that $v \in \big \langle \ker x_0^\top \: | \: X \big \rangle$. 
	
Conversely, suppose that $v \in \big \langle \ker x_0^\top \: | \: X \big \rangle$. This implies that $x_0^\top e^{Xt} v = 0$ for all $t \in \mathbb{R}_{\geq 0}$. We compute
\begin{equation}
\begin{aligned}
P v &= \int_0^{T} e^{Xt} x_0 x_0^\top e^{Xt} v \: dt = 0.
\end{aligned}
\end{equation}
In other words, we obtain $v \in \ker P$. We conclude that $\ker P = \big \langle \ker x_0^\top \: | \: X \big \rangle$, which completes the proof.
\end{IEEEproof}

\begin{IEEEproof}[Proof of Theorem \ref{theorem unique solution}]
	To prove the `if' part, suppose that the network reconstruction problem is not solvable for $(x_0,X,\mathcal{K})$. We want to prove that the solution to \eqref{lyapunov equation theorem1} is not unique. By hypothesis, there exists a matrix $\bar{X} \in \im \mathcal{K} \setminus \{X\}$ such that $e^{Xt}x_0 = e^{\bar{X}t} x_0$ for all $t \in [0,T]$. We can repeat the discussion of Equation \eqref{derivation} for $\bar{X}$, to show that $\bar{X}$ also solves the Lyapunov equation \eqref{lyapunov equation theorem1}. Consequently, we conclude that there exists no unique solution $S$ satisfying \eqref{lyapunov equation theorem1}.
	
	Conversely, to prove the `only if' part, suppose that there exists no unique solution to \eqref{lyapunov equation theorem1}. Note that $S = X$ is always a solution to \eqref{lyapunov equation theorem1} by Equation \eqref{lyapunov equation transformed}. This implies that there exists a matrix $\bar{X} \neq X$ satisfying \eqref{lyapunov equation theorem1}. Consequently, $\bar{X} \in \im \mathcal{K}$, and
	\begin{equation}
	\label{differenceLyapunov}
	(\bar{X} - X)P + P(\bar{X} - X) = 0.
	\end{equation}
	Since $P$ is symmetric positive semidefinite, there exists an orthogonal matrix $U \in \mathbb{R}^{n \times n}$ such that $P = U \Lambda U^\top$, where 
	\begin{equation*}
	\Lambda = \begin{pmatrix}
	D & 0 \\ 0 & 0 
	\end{pmatrix},
	\end{equation*}
	with $D$ a positive definite diagonal matrix. We define the matrix $\hat{X} := U^\top (\bar{X} - X)U$. It follows from \eqref{differenceLyapunov} that $\hat{X}$ satisfies the Lyapunov equation $\hat{X}\Lambda + \Lambda \hat{X} = 0$. Next, we partition $\hat{X}$ as 
	\begin{equation*}
	\hat{X} = \begin{pmatrix}
	\hat{X}_{11} & \hat{X}_{12} \\ \hat{X}_{21} & \hat{X}_{22} 
	\end{pmatrix},
	\end{equation*}
	where the partitioning of $\hat{X}$ is compatible with the one of $\Lambda$. Then, we rewrite $\hat{X}\Lambda + \Lambda \hat{X} = 0$ as
	\begin{equation*}
	\begin{pmatrix}
	\hat{X}_{11} D + D \hat{X}_{11} & D \hat{X}_{12} \\
	\hat{X}_{21} D & 0
	\end{pmatrix} = \begin{pmatrix}
	0 & 0 \\ 0 & 0
	\end{pmatrix}.
	\end{equation*}
	Since $D$ is nonsingular, $\hat{X}_{12} = 0$. Moreover, since $D$ and $-D$ do not have common eigenvalues, the Lyapunov equation $\hat{X}_{11} D + D \hat{X}_{11} = 0$ has a unique solution given by $\hat{X}_{11} = 0$ (cf. Theorem 2.5.10 of \cite{Basile1992}). This means that $\Lambda \hat{X} = 0$. Therefore, $P(\bar{X} - X) = 0$. By Proposition \ref{proposition kernel} we have $x_0^\top X^i (\bar{X} - X) = 0$ for $i = 0,1,\dots,n-1$. By exploiting symmetry, we obtain $x_0 \in \big \langle \ker(\bar{X} - X) \mid X \big \rangle$. We conclude by Theorem \ref{theorem nec suf} that the network reconstruction problem is not solvable for $(x_0,X,\mathcal{K})$. 
	
	Finally, as we have shown in Equation \eqref{lyapunov equation transformed} that $X \in \im \mathcal{K}$ is always a solution to the Lyapunov equation $SP + PS = Q$, it is immediate that $S = X$ if there exists a unique solution $S$ to \eqref{lyapunov equation theorem1}.
\end{IEEEproof}

Theorem \ref{theorem unique solution} provides a general framework for network reconstruction. Indeed, suppose that the network reconstruction problem is solvable for $(x_0,X,\mathcal{K})$. We can compute the matrices $P$ and $Q$ from the state measurements $x_{x_0}(t)$ for $t \in [0,T]$. Then, network reconstruction boils down to computing the unique solution $S$ to the constrained Lyapunov equation \eqref{lyapunov equation theorem1}. In the subsequent sections, we will show how this can be done for several types of network dynamics.

\subsection{Network reconstruction for $\mathcal{K} = \mathcal{Q}$}
\label{section identification Q}

In this section, we consider network reconstruction in the case that $\mathcal{K}$ is equal to $\mathcal{Q}$. Based on Theorem \ref{theorem unique solution} we will derive an algorithm to identify the unknown matrix $X \in \mathcal{Q}(G)$ using state measurements taken from the network. 

Recall from Theorem \ref{theorem unique solution} that the network reconstruction problem is solvable for $(x_0,X,\mathcal{Q})$ if and only if there exists a unique matrix $S$ satisfying $SP + PS = Q$ and $S \in \im \mathcal{Q}$.
Note that $\im \mathcal{Q}$ is equal to $\mathbb{S}^n$, the set of $n \times n$ symmetric matrices. In other words, if the network reconstruction problem is solvable for $(x_0,X,\mathcal{Q})$, the solution to the problem can be found by computing the unique symmetric solution to the Lyapunov equation $SP + PS = Q$. It is not difficult to see that there exists a unique \emph{symmetric solution} to $SP + PS = Q$ if and only if there exists a unique \emph{solution} to $SP + PS = Q$. This yields the following corollary of Theorem \ref{theorem unique solution}.

\begin{corollary}
\label{corollary Q}
Let $G \in \mathcal{G}^n$ be a graph, and let $X \in \mathcal{Q}(G)$. Moreover, consider a vector $x_0 \in \mathbb{R}^{n}$, and let $P$ and $Q$ be as defined in \eqref{transformed P} and \eqref{lyapunov equation transformed} respectively. The network reconstruction problem is solvable for $(x_0,X,\mathcal{Q})$ if and only if the Lyapunov equation $SP + PS = Q$ admits a unique solution $S$. Under this condition, we have $S = X$.
\end{corollary}

Based on Corollary \ref{corollary Q}, we establish Algorithm 1, which infers the state matrix $X$ and graph $G$ from measurements. Recall from Theorem \ref{theorem Q nec suf} that the network reconstruction problem is solvable for $(x_0,X,\mathcal{Q})$ if and only if $(x_0^\top,X)$ is observable. Of course, we can not directly check observability of $(x_0^\top,X)$ since $X$ is unknown. However, we can in fact check observability of the pair $(x_0^\top,X)$ using the matrix $P$. Indeed, by Proposition \ref{proposition kernel}, $(x_0^\top,X)$ is observable if and only if the matrix $P$ is nonsingular. Note that this condition is similar to the so-called \emph{persistency of excitation} condition, found in the literature on adaptive systems (cf. Section 3.4.3 of \cite{Mareels1996}).
\begin{algorithm}[H]
	\caption{Network reconstruction for $(x_0,X,\mathcal{Q})$}
	\begin{algorithmic}[1]
		\Require{Measurements $x_{x_0}(t)$ for $t \in [0,T]$};
		\Ensure{Matrix $X$ or ``No unique solution exists"};
		\State Compute the matrix $P = \int_0^{T} x_{x_0}(t) x_{x_0}(t)^\top dt$;
		\If{$\rank P < n$}
		\State \Return{``No unique solution exists"};
		\Else  
		\State Compute the matrix $Q = x_0 x_0^\top - x_T x_T^\top$;
		\State Solve $SP + PS = Q$ with respect to $S$;
		\State \Return{$X = S$};
		\EndIf
	\end{algorithmic}
\end{algorithm}
A classic method to solve the Lyapunov equation in Step 6 of Algorithm 1 is the \emph{Bartels-Stewart algorithm} \cite{Bartels1972}. In addition, much effort has been made to develop methods for solving large-scale Lyapunov equations \cite{Haber2016}, \cite{Simoncini2007}. Typically, such methods use the Galerkin projection of the Lyapunov equation onto a lower-dimensional Krylov subspace \cite{Simoncini2007}. The resulting reduced problem is then solved by means of standard schemes for (small) Lyapunov equations. Using these techniques, it is possible to efficiently solve large-scale ($n > 10000$) Lyapunov equations \cite{Simoncini2007}.

\begin{remark}
In theory, the correctness of Theorem \ref{theorem unique solution}, Corollary \ref{corollary Q}, and Algorithm 1 is independent of the exact choice of time $T > 0$. However, choosing small $T$ results in a matrix $P$ with high condition number, and hence numerical rank computation (as in line 2 of Algorithm 1) becomes inaccurate. Consequently, in practice the value of $T$ should be sufficiently large. In our simulations, good results were obtained using $T = 10$ (see Section \ref{section illustrative example}).
\end{remark}

\begin{remark}
Even though the focus of this note is on continuous-time systems, we remark that Algorithm 1 can also be applied for network reconstruction of \emph{discrete-time} networks of the form 
\begin{equation}
\label{discrete time}
\begin{aligned}
z(k+1) &= Mz(k) \text{ for } k \in \mathbb{N} \\
z(0) &= z_0,
\end{aligned}
\end{equation}
where $z \in \mathbb{R}^n$ and $M \in \im \mathcal{Q}$. In this case, we assume that we can measure the state $z(k)$, for $k = 0,1,\dots,m$, where $m \geq n$. From these measurements, we compute
\begin{equation*}
P := \sum_{k=0}^{m-1} z(k)z(k)^\top, \:\:\: Q := \sum_{k=0}^{m-1} z(k+1)z(k)^\top+z(k)z(k+1)^\top.
\end{equation*}
Similar to the continuous-time case, the matrix $P$ is nonsingular if and only if $(z_0^\top,M)$ is observable. Under this condition, we can reconstruct $M$ by computing the unique solution to the Lyapunov equation $MP + PM = Q$.

The above approach can also be used for the continuous-time network \eqref{networked system} in the case that we cannot measure the state trajectory $x_{x_0}(\D)$ during a \emph{time interval}, but only have access to \emph{sampled measurements}. Indeed, suppose that we can measure $x_{x_0}(k\tau)$ for $k = 0,1,\dots,m$, where $\tau \in \mathbb{R}_{>0}$ is some sampling period. We can then use the framework for discrete-time systems on $z(k) := x_{x_0}(k\tau)$ to reconstruct the matrix $M = e^{X\tau}$. Subsequently, we can reconstruct $X$ by computing the (unique) matrix logarithm of $e^{X\tau}$.
\end{remark}

\subsection{Network reconstruction for $\mathcal{K} = -\mathcal{L}$ and $\mathcal{K} = \mathcal{A}$}
\label{section identification L}
Although Algorithm 1 is applicable to general network dynamics described by state matrices $X \in \mathcal{Q}(G)$, the observability condition guaranteeing uniqueness of the solution to \eqref{lyapunov equation transformed2} can be quite restrictive if the \emph{type} of network is a priori known. We have already seen in Section \ref{section solvability L} that observability of the pair $(x_0^\top,X)$ is not necessary for the solvability of the network reconstruction problem for adjacency or consensus networks. Therefore, in this section we focus on network reconstruction for $(x_0,-L,-\mathcal{L})$ and $(x_0,A,\mathcal{A})$. 

Recall from Theorem \ref{theorem unique solution} that the network reconstruction problem is solvable for $(x_0,-L,-\mathcal{L})$ if and only if there exists a unique matrix $S$ satisfying $SP + PS = Q$ and $S \in -\im \mathcal{L}$. Based on the definition of $\mathcal{L}$ (see Section \ref{preliminaries graphs}), we find the following corollary of Theorem \ref{theorem unique solution}.
\begin{corollary}
	\label{corollary L}
	Let $G \in \mathcal{G}^n$ be a graph, and let $L \in \mathcal{L}(G)$. Moreover, consider a vector $x_0 \in \mathbb{R}^{n}$, and let $P$ and $Q$ be as defined in \eqref{transformed P} and \eqref{lyapunov equation transformed} respectively. The network reconstruction problem is solvable for $(x_0,-L,-\mathcal{L})$ if and only if there exists a unique solution $S$ to
	\begin{equation}
	\label{constrained lyapunov L}
	\begin{aligned}
	SP + PS &= Q, \quad S \in \mathbb{S}^n, \quad S \mathbbm{1} = 0, \quad S_{ij} \geq 0 \text{ for } i \neq j. 
	\end{aligned}
	\end{equation}  	
Moreover, under this condition, we have $S = -L$.
\end{corollary}
The constraint $S_{ij} \geq 0$ for $i \neq j$ can be stated as a \emph{linear matrix inequality} (LMI) in the matrix variable $S$. Indeed, $S_{ij} \geq 0$ is equivalent to $e_i^\top S e_j \geq 0$, where $e_k$ denotes the $k$-th column of the $n \times n$ identity matrix. Consequently, by Corollary \ref{corollary L}, network reconstruction for $(x_0,-L,-\mathcal{L})$ boils down to finding the matrix $S$ satisfying linear matrix equations and linear matrix inequalities, given by \eqref{constrained lyapunov L}. There is efficient software available to solve such problems. See, for instance, the LMI Lab package in Matlab and Yalmip \cite{Lofberg2005}. We can deduce a corollary similar to Corollary \ref{corollary L} for the class $\mathcal{A}(G)$. In this case, the restrictions on the elements of $S$ are $S_{ii} = 0$ and $S_{ij} \geq 0$ for all $i \in V$ and all $j \neq i$.  

\section{Illustrative example}
\label{section illustrative example}
\begin{figure}[h!]
	\includegraphics[width=0.47\textwidth]{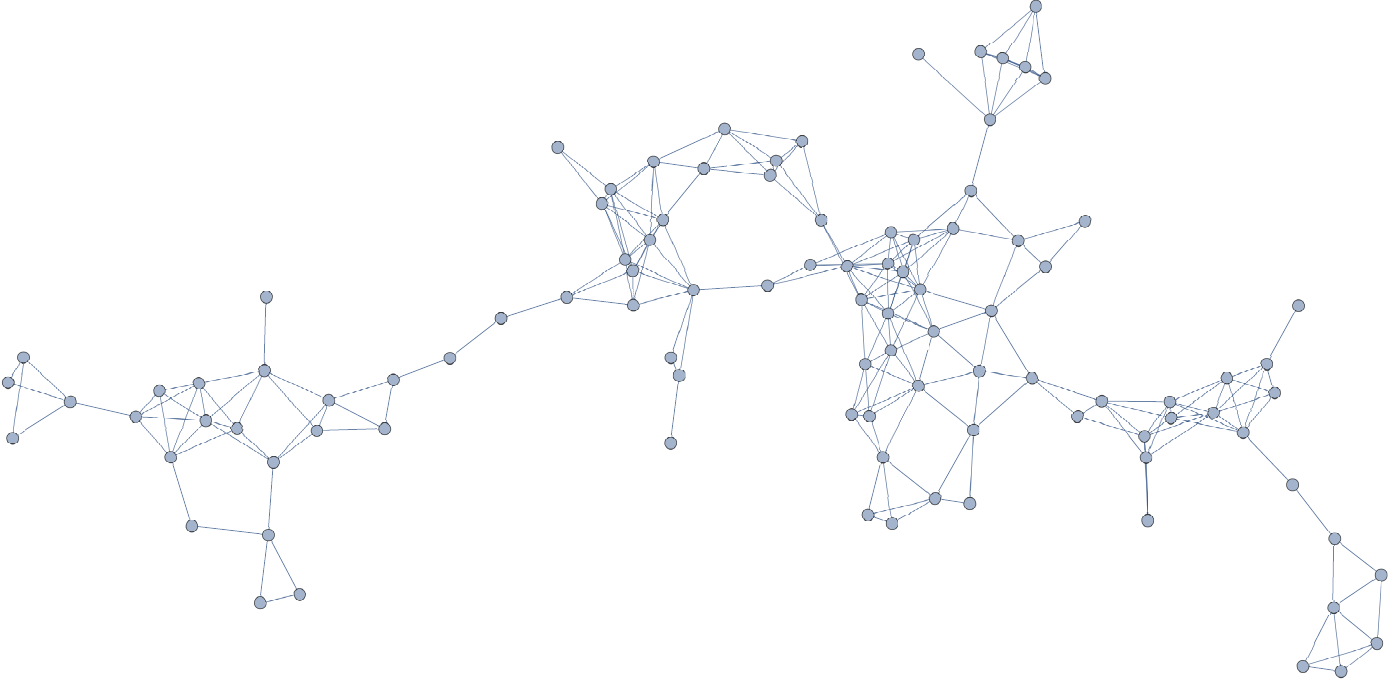}
	\caption{Graph $G$ of the sensor network.}
	\label{fig:sensornetwork}
\end{figure}

In this section we illustrate the developed theory by consi\-dering an example of a sensor network. Specifically, consider a graph $G = (V,E)$ consisting of 100 sensor nodes, monitoring a region of $1 \: \mathrm{km} \times 1 \: \mathrm{km}$ (see Figure \ref{fig:sensornetwork}). It is assumed that the sensors are linked using a so-called \emph{geometric link model} \cite{Penrose2003}. This means that there is a connection between two nodes in the network if and only if the distance between the two nodes is less than a certain threshold, set to be equal to $135 \: \mathrm{m}$ in this example. It is assumed that the sensors run consensus dynamics, that is, the dynamics of the network is given by $\dot{x}(t) = -Lx(t)$, where $x \in \mathbb{R}^{100}$, and $L \in \mathcal{L}(G)$ is the unweighted Laplacian associated with $G$. The components of the initial condition $x_0 \in \mathbb{R}^{100}$ were selected randomly within $[0,10]$. Moreover, for this example, measurements were used over the time-interval $[0,10]$, i.e., $T = 10$. We compute the matrices $P$ and $Q$, and solve \eqref{constrained lyapunov L} using Yalmip. The resulting identified Laplacian matrix is denoted by $L_r$. The relative and maximum element-wise errors between the identified Laplacian $L_r$ and original Laplacian $L$ are very small. Specifically, we obtain
\begin{equation*}
\frac{\norm{L_{r}-L}}{\norm{L_{r}}} = 1.56 \cdot 10^{-8}, \quad \max_{i,j \in V} |L_{ij}-(L_r)_{ij}| = 2.21 \cdot 10^{-7},
\end{equation*}
where $\| \D \|$ denotes the induced 2-norm.

\section{Conclusions}
\label{section conclusions}

In this technical note, we have considered the problem of network reconstruction for networks of linear dynamical systems. In contrast to papers studying network reconstruction for specific network dynamics such as consensus dynamics \cite{Nabi-Abdolyousefi2010} and adjacency dynamics \cite{Fazlyab2014}, we considered network reconstruction for general linear network dynamics described by state matrices contained in the qualitative class. 
We formulated what is meant by solvability of the network reconstruction problem. Subsequently, we provided necessary and sufficient conditions under which the network reconstruction problem is solvable. Using constrained Lyapunov equations, we established a general framework for network reconstruction of networks of dynamical systems. We have shown that this framework can be used for a variety of network types, including consensus and adjacency networks. Finally, we have illustrated the theory by reconstructing the network topology of a sensor network. 

\ifCLASSOPTIONcaptionsoff
 \newpage
\fi

\bibliographystyle{IEEEtran}
\bibliography{MyRef}

\end{document}